\documentstyle[12pt]{article}
\def\d{\partial}
\def\t{\theta}
\def\a{\alpha}
\def\P{\Pi}
\def\lb{\label}
\def\be{\begin{equation}}
\def\ee{\end{equation}}
\begin{document}

\title{ Paragrassmann Integral, Discrete Systems
and Quantum Groups.}
\author{A.P.Isaev \\
{\it Bogoliubov Laboratory of Theoretical Physics,} \\
{\it JINR, Dubna, Moscow Region, 141980, Russia}}
\date{}
\maketitle

\vspace{0.5cm}

\begin{center}
{\bf Abstract.}
\end{center}

This report is based on review paper \cite{FIIK}.
Some aspects of differential and integral calculi
on generalized grassmann (paragrassmann) algebras are considered.
The integration over paragrassmann variables is
applied to evaluate the partition function
for the $Z_{p+1}$ Potts model on a chain.
Finite dimensional paragrassmann representations
for $GL_{q}(2)$ are constructed.

\vspace{1cm}

Generalizations of grassmann algebras have been considered
by many authors (see e.g. \cite{PA}-\cite{PA2}
and references therein)
>from different points of view.
Those generalizations were stimulated by
investigations in 2D conformal field theories
\cite{PA1,Zam}, anionic models and topological
field theories which led to consideration of
unusual statistics. The last one includes not only well known
parastatistics \cite{PARA} but also fractional and braid statistics
\cite{PA2,LM}. We note also attempts to generalize
supersymmetry to parasupersymmetry \cite{MA,PA1}.
Our construction of paragrassmann
algebras (PGA) with many variables
\cite{fik,fik2} is a direct generalization
of the Weyl construction \cite{W} for the grassmann
Heisenberg-Weyl algebra.

Let us consider the algebra $\Pi_{p+1}$ with two nilpotent
generators $\t$ and $\d$:
\be
\lb{1.1}
\t^{p+1} = 0 = \d^{p+1} \; , \;\; \t^{p} \neq 0 \; , \;\;
\d^{p} \neq 0 \; ,
\ee
where $p$ is a positive integer
(usually called the order of parastatistics).
The defining relation for this algebra is chosen
such that one can push $\d$ to the right:
\be
\lb{1.2}
\d \t = b_{0} + b_{1} \t \d + b_{2} \t^{2} \d^{2} + \dots +
b_{p} \t^{p} \d^{p} \; .
\ee
Here $b_{i}, \; b_{0} \neq 0$ are complex numbers restricted by
consistency of eqs. (\ref{1.1}) and (\ref{1.2}). For
general values of $\{ b_{i} \}$ one can transform the
derivative as
$$
\d \rightarrow \sum_{n=1}^{p} c_{n} \t^{n-1} \d^{n} \; ,
$$
($c_{i}$ are parameters)
and reduce relations (\ref{1.1}), (\ref{1.2})
to the simplest form
($b_{0}=1$, $b_{1}=q$ and $b_{k} = 0$ for $k \geq 2$)
\be
\lb{1.8}
\d \t = 1 + q \t \d \; , \;\; \d^{p+1} = 0 = \t^{p+1} \; ,
\ee
We call this version of PGA the $q$-version.
In this case, we deduce
\be
\lb{1.8a}
\d \, \t^{i} = \alpha_{i} \, \t^{i-1} + q^{i} \, \t^{i} \, \d  \; ,
\ee
where
$$
\a_{i}=1+q+ \dots +q^{i-1} = \frac{1-q^{i}}{1-q} \equiv (i)_q
\equiv i_{q}\; .
$$
Equation $\d \t^{p+1} =0$ is equivalent to the condition
$\alpha_{p+1}=0$ or $q^{p+1}=1  \;\; (q \neq 1)$
while the condition $\d^{p} \neq 0$ requires $\a_{i} \neq 0$ or
$q^{n+1} \neq 1$ $(n<p)$, i.e.,
$q$ must be a primitive root of unity.
Note that for $p=1$
algebra (\ref{1.8}) is isomorphic to the grassmann Heisenberg-Weyl
algebra.

The well known relations (\ref{1.8}) and (\ref{1.8a})
have been introduced in \cite{fik} by assuming
that $\d$ is a differential operator satisfying
a generalized Leibniz rule.
It is clear that the $q$-version of
$\P_{p+1}$ (\ref{1.8}) is related to
$q$-oscillators and quantum groups
(see \cite{fik} and references therein).
This version is also useful for constructing
PGA's with many generators $\t_{j}$ and $\d_{j}$
(see below). The other useful versions of PGA
(\ref{1.1}), (\ref{1.2}) have been discussed
in \cite{fik2}.

{} Fundamental representations of the operators
$\t$ and $\d$ in the $q$-version have the form:
\be
\lb{1.9a}
\t_{mn}  =  \langle m| \t |n\rangle =\beta_{n+1} \delta_{m,n+1}
\; \; , \;\;\;
\d_{mn}  =  \langle m| \d |n\rangle =
\frac{n_q }{ \beta_{n} }\delta_{m,n-1} \; \; .
\ee
where $m,n =0,1,\dots ,p$; $\beta_i \neq 0$ are
arbitrary parameters, and we introduce
the ladder of $p+1$ states $|n>$ defined by
$$
\d |0\rangle = 0 \;, \; |k\rangle \sim \t^{k} |0\rangle \; ,
\t |k\rangle = \beta_{k+1} |k+1\rangle \; .
$$
The dual states $<m|$ satisfy the conditions
$<m|n> = \delta_{m,n}$. From definition of vacuum $|0>$
one can obtain relations
\be
\lb{1.19}
< 0 | \d^{n} \, \t^{m} | 0 > =
\delta_{n,m} \, n_q !   \; ,
\ee
which will be needed below. We introduce also the operator
\be
\lb{1.10a}
g=\d \, \t - \t \, \d \;\; \Longrightarrow \;\; g^{p+1}=1 \; ,
\ee
which generates automorphisms in PGA (\ref{1.8}):
$$
g \, \t \, g^{-1} = q \, \t  \; , \;\;
g \, \d \, g^{-1} = q^{-1} \, \d  \; .
$$
Using operator $g$ we define the natural hermitian conjugation
for the elements $\t$ and $\d$ of PGA:
\be
\lb{1.13}
\t^{\dagger} =  C \, g^{-1/2} \, \d \;\; , \;\;\;
\d^{\dagger} = \t \, g^{-1/2} \, C^{-1} \; ,
\ee
where $g^{1/2} \, \t = q^{1/2} \, \t \, g^{1/2}$,
$g^{1/2} \, \d = q^{-1/2} \, \d \, g^{1/2}$ and invertible element
$C \in \P_{p+1}$ is chosen such that
$[C, \, \t \, \d ] =0$. From eqs. (\ref{1.10a}),
(\ref{1.13}) one can obtain
$g^{\dagger} = g^{-1}$, and the conditions
$\t^{\dagger \dagger} = \t$, $\d^{\dagger \dagger} = \d$
yield the relation $C = C^{\dagger}$.
Now the first relation in (\ref{1.8}) is rewritten in
the form
\be
\lb{1.14}
\t^{*} \, \t -  q^{1/2} \, \t \, \t^{*} = g^{-1/2}  \; , \;\;
\t^{*} \equiv C^{-1} \, \t^{\dagger} \; .
\ee
As it has been shown in \cite{FIIK}, one can
put $C=1$ for $q=\exp (2 \pi i/(p+1))$
and relations (\ref{1.14}) coincide with the definition
of Macfarlane-Biedenharn's $q$-oscillator \cite{MaBi}. For
other roots of unity
$q=\exp (2 \pi i n/(p+1))$ $(n=2,\dots,p)$,
the condition of positive definiteness
$<0|(\t^{\dagger})^{n} \, \t^{n} |0> \, \geq \, 0$
leads to the statement $C \neq 1$ \cite{FIIK}. In this case
relations (\ref{1.14}) define the modification
of Macfarlane-Biedenharn's $q$-oscillator.

Now we introduce an integral over paragrassmann (PG) variables
which generalizes Berezin's integral over grassmann variables.
Let us consider the algebra of
PG numbers $\t$ and $\overline{\t}$:
\be
\lb{1.20}
\overline{\t} \, \t = q \, \t \, \overline{\t} \; , \;\;
\overline{\t}^{p+1} = 0 = \t^{p+1} \; .
\ee
An element of this algebra is represented as a
function of $\t$ and $\overline{\t}$
\be
\lb{1.21a}
f(\t, \, \overline{\t}) =
\sum_{n,m=0}^{p} \, \frac{1}{n_{q}!} \, f_{nm}\,
\t^{n} \, \overline{\t}^{m} \; ,
\ee
where $f_{ij}$ is a matrix of complex numbers.

{} Functions which depend only on
$\t$ (or only on $\overline{\t}$) are called holomorphic
\be
\lb{1.22}
g(\t) =
\sum_{n=0}^{p} \, g_{n} \, \t^{n} \equiv <\t | g> \;\; , \;\;\;
f(\overline{\t}) =
\sum_{n=0}^{p} \, f_{n} \, \overline{\t}^{n} \equiv
<f | \overline{\t}>\; .
\ee
The natural definition of the left derivative with respect to
$\t$ can be deduced from (\ref{1.8a}) and looks like
$$
\d \triangleright (g(\t)) =
\sum_{n=1}^{p} \, (n)_q \, g_{n} \, \t^{n-1} \; .
$$
The requirement
$\int d\t \, \d \triangleright (g(\t))  = 0$
gives us the generalization of Berezin's integral for the
PG case
$$
\int d\t \, \t^{n}  = x_{p} \, \delta_{n,p} \; ,
$$
where $x_{p}$ is a complex number which will be fixed
below. Analogous integral is postulated for
variables $\overline{\t}$:
$$
\int d\overline{\t}
\, \overline{\t}^{n} =
\overline{x}_{p} \, \delta_{n,p} \Leftrightarrow
\int  d\overline{\t}
\, \overline{\d} \triangleright (f(\overline{\t}))
= 0 \; .
$$
Our aim is to find integral representation for the
matrix elements (\ref{1.19}), i.e.
\be
\lb{1.25}
<0| \, \d^{n} \, \t^{m} \, |0> =
\int \int d\t \, d\overline{\t}
\; \overline{\t}^{n} \, \t^{m} \mu (\t, \, \overline{\t})
= \delta_{n,m} \, n_q ! \; ,
\ee
where $\mu (\t, \, \overline{\t})$ is a measure function
which must be defined.
The 2-fold integral in (\ref{1.25})
is understood as iterated:
$$
\int d\t \left( \int d \overline{\t} f(\overline{\t}) \right) g(\t) \; .
$$

Note that for the primitive root $q$ of unity
we have the following identity:
$$
p_q ! = (-1)^{n} \, q^{-n(n+1)/2} \, n_q ! \, (p-n)_{q} ! \; .
$$
Using this identity one can find that relation
(\ref{1.25}) is fulfilled only
if we restrict numbers $x_{p}$ and $\overline{x}_{p}$
by the relation
$x_{p} \, \overline{x}_{p} = (p)_q !$
and choose the function $\mu$ in the form
\be
\lb{1.28}
\mu(\t, \, \overline{\t}) =
\sum_{n=0}^{p} \frac{(-\t \, \overline{\t})^{n}}{(n)_q !}  \equiv
exp_{q} (-\t \, \overline{\t}) \; .
\ee
Thus, for two arbitrary holomorphic functions
$f(\d)$ and $g(\t)$ we have the following integral
representation:
\be
\lb{1.29}
 <0| \, f(\d) \, g(\t) \, |0> =
\int \int d\t \, d\overline{\t}
\, f(\overline{\t}) \, g(\t) \, exp_{q}(- \t \, \overline{\t}) \,
\; .
\ee
Then, from the natural conditions
\be
\lb{1.29c}
\int \int d\t \, d\overline{\t}
\, f(\overline{\t}) \, \d \triangleright \left( g(\t) \, \mu
( \t \, \overline{\t}) \right) \,
= 0 =
\int \int d\t \, d\overline{\t}
\, \overline{\d} \triangleright \left( f(\overline{\t}) \,
g(\t) \, \mu ( \t \, \overline{\t}) \right) \; ,
\ee
we derive (see \cite{FIIK}) the whole set of defining relations for the
differential algebra $\Pi_{p+1}(2)$ on the PG plane (\ref{1.20})
\be
\lb{1.29a}
\d \, \overline{\t} = q \, \overline{\t} \, \d \;  , \;\;
\overline{\d} \, \overline{\t} = 1 + q \, \overline{\t} \, \overline{\d}
\; , \;\;
\t \, \overline{\d} = q \, \overline{\d} \, \t \; ,  \;\;
\overline{\d} \, \d = q \, \d \, \overline{\d} \; .
\ee
Here, we have to add the relations
\be
\lb{a1.8}
\d \t = 1 + q \t \d \; , \;\; \d^{p+1} = 0 = \overline{\d}^{p+1} \; ,
\ee
postulated earlier.

Note that integrals of the type (\ref{1.25}) and (\ref{1.29})
for the nilpotent quantum plane (\ref{1.20})
have also been considered in \cite{BF}, where
a slightly different differential algebra $\P_{p+1}(2)$
was proposed.

Taking into account relations (\ref{1.29a}), (\ref{a1.8})
we remark that evaluation of
integral (\ref{1.29}) is equivalent to the action of
a differential operator $(\d^{p} \, \overline{\d}^{p})$
\be
\lb{1.29e}
\int \int d\t \, d\overline{\t}
\, f(\overline{\t}) \, g(\t) \, exp_{q}(- \t \, \overline{\t}) \,
=
\frac{1}{p_{q}!} (\d^{p} \, \overline{\d}^{p}) \triangleright
\left( f(\overline{\t}) \, g(\t)
\, \mu(\t \overline{\t}) \right) \; .
\ee
The next point is that the differential algebra
$\Pi_{p+1}(2)$ on the PG plane
can be constructed as a direct product of
two differential algebras $\Pi_{p+1}(1)$ (\ref{1.8}) on the PG line.
Indeed, one can check that the operators
\be
\lb{1.36}
\t_{1} = \t \otimes 1  \; , \;\;
\overline{\t}_{1} = g \otimes \t \; , \;\;
\d_{1} = \d \otimes 1 \; ,
\overline{\d}_{1} = g^{-1} \otimes \d \; ,
\ee
satisfy eqs. (\ref{1.20}), (\ref{1.29a})
and (\ref{a1.8}).
The explicit form of representations (\ref{1.36})
gives us the idea how
to extend the algebra $\Pi_{p+1}(2)$
up to the PG algebra $\Pi_{p+1}(2N)$ with the $2N+2N$ generators
$\{ \t_{i}, \; \d_{i}, \; \overline{\t}_{i}, \; \overline{\d}_{i} \}$
where $i=1, \dots , N$.
Indeed, by analogy with formulas (\ref{1.36}),
generators of $\Pi_{p+1}(2N)$ can be realized as:
\be
\lb{1.35}
\begin{array}{c}
\t_{i} = G^{a_{1}} \otimes \dots \otimes
G^{a_{i-1}} \otimes (\t \otimes 1)
\otimes I^{\otimes (N-i)} \; , \\
\overline{\t}_{i} = G^{b_{1}} \otimes \dots \otimes
G^{b_{i-1}} \otimes (g \otimes \t)
\otimes I^{\otimes (N-i)} \; , \\
\d_{i} = G^{-a_{1}} \otimes \dots \otimes G^{-a_{i-1}} \otimes (\d \otimes 1)
\otimes I^{\otimes (N-i)} \; , \\
\overline{\d}_{i} = G^{-b_{1}} \otimes \dots \otimes
G^{-b_{i-1}} \otimes (g^{-1} \otimes \d)
\otimes I^{\otimes (N-i)} \; ,
\end{array}
\ee
where
$I = 1 \otimes 1$,
$G^{a_{i}}=g^{a^{1}_{i}} \otimes g^{a^{2}_{i}}$, and
components of the 2-D vectors $a_{i}$, $b_{i}$ are equal to
$a^{\alpha}_{j}, \, b^{\alpha}_{k} = \pm 1$.

Now, we show how one can apply PGA
for the consideration of the
$Z_{p+1}$ Potts model on a closed 1-dimensional
chain with $N$ cites. This model is formulated in the following
manner \cite{Bax}. We assign spin $\sigma_{i}$
to each cite $i$ of the chain. This spin takes one of
the $p+1$ values, say $0, 1, 2, ..., p$. We consider
the closed chain and identify $\sigma_{N+1} = \sigma_{1}$.
Then, the partition function for this model is defined as
\be
\lb{3.2}
Z(N) = \sum_{\sigma} \exp \{ K \sum_{n=1}^{N}
\delta(\sigma_{n},\sigma_{n+1})  \} \; .
\ee
Here
\be
\lb{3.1}
\delta(\sigma_{i},\sigma_{j})
 = 1 \;\; for \;\; \sigma_{i}=\sigma_{j} \;\; and \;\;
\delta(\sigma_{i},\sigma_{j})
 = 0 \;\; for \;\; \sigma_{i} \neq \sigma_{j} \; ,
\ee
$K = J/T$, $J$ is a constant and
$T$ - temperature.
The partition function (\ref{3.2})
can be rewritten in terms of transfer-matrices (see \cite{Bax}):
\be
\lb{3.3}
Z(N) = Tr(V^{N}) =
\sum_{\sigma}
V_{\sigma_{1},\sigma_{2}} \, V_{\sigma_{2},\sigma_{3}} \dots
V_{\sigma_{N},\sigma_{1}} \; ,
\ee
where a transfer-matrix $V$ has the form
\be
\lb{3.4}
V_{\sigma_{n},\sigma_{n+1}} =
\left( \, 1 +  (x-1) \delta(\sigma_{n},\sigma_{n+1}) \, \right)
\; , \;\; (x=\exp(K)) \; .
\ee
Our aim is to demonstrate that the partition function
(\ref{3.2}), (\ref{3.3}) can be calculated
by representing $Z(N)$ in the form
of $2N$-fold integrals over PG variables
introduced above. As by product we obtain additional
conditions on the defining relations of the algebra
$\Pi_{p+1}(2N)$, which fix this algebra completely.
To perform this, we use the idea of grassmann
factorization method successfully applied
to the solution of the 2-dimensional Ising
model (see \cite{P}). More correctly, we formulate the direct
generalization of this method to the case of PG variables.

{} First of all, we note that the
$\delta$- function (\ref{3.1}) is represented
as a finite sum
$$
\delta(\sigma, \, \sigma') =
\frac{1}{p+1} \, \sum_{m=0}^{p} \, q^{m(\sigma - \sigma')} \; ,
$$
where $q= \exp(2 \, i \, \pi / (p+1))$ - a primitive root of
unity. After this, the transfer-matrix
(\ref{3.4}) is written for the $n$-th cite
in the form
\be
\lb{3.10}
V_{\sigma_{n},\sigma_{n+1}} =
\int \int d\t_{n} d\overline{\t}_{n} \, \left[ \mu(\t_{n} \overline{\t}_{n})
(\sum_{m=0}^{p} \frac{t_{m}}{m_{q}!} q^{m\sigma_{n}} \overline{\t}_{n}^{m} )
(\sum_{k=0}^{p} q^{-k \sigma_{n+1}} \t_{n}^{k} ) \right] \; ,
\ee
where we use formula (\ref{1.25}) and introduce the constants
$$
t_{0} = \frac{(p+x)}{(p+1)} \; , \;\;
t_n = \frac{(x-1)}{(p+1)} \;\; (n \geq 1) \; .
$$
Now let us substitute (\ref{3.10}) in (\ref{3.2})
and perform summation over spins $\sigma_{i}$ with the help of identity
($0 \leq k \leq p$)
$$
\sum_{\sigma=0}^{p} \, q^{k\sigma} = (p+1) \, \delta_{k,0} \; .
$$
As a result, we obtain the product of
$N$ $2$-fold integrals (\ref{3.10}) which
we would like to rewrite into one $2N$-fold
integral over PG variables. It can be done if
we choose the vectors $a_{i} = - b_{i}$ in the definition
of the algebra $\Pi_{p+1}(2N)$ (\ref{1.35}).
In this case, the partition function (\ref{3.3}) takes the form
$$
Z(N) = (p+1)^{N} \, \int \dots \int d\t_{1} d\overline{\t}_{1} \dots
d\t_{N} d\overline{\t}_{N} \, \exp_{q}(-\t_{1}\overline{\t}_{1})
 \cdots \exp_{q}(-\t_{N} \, \overline{\t}_{N}) \,
$$
\be
\lb{3.12}
\sum_{m_{N}} \frac{t_{m_{N}}}{(m_{N})_{q}!}
\overline{\t}_{1}^{m_{N}} \, \left(
\sum_{m_{1}} \frac{t_{m_{1}}}{(m_{1})_{q}!}
\t_{1}^{m_{1}}\overline{\t}^{m_{1}}_{2}
\, \dots \,
\sum_{m_{N-1}} \frac{t_{m_{N-1}}}{(m_{N-1})_{q}!}
\t_{N-1}^{m_{N-1}}\overline{\t}^{m_{N-1}}_{N} \right)
\,
\t^{m_{N}}_{N} \; .
\ee
Using formulas (\ref{1.25}) and
(\ref{1.29e}) one can calculate all PG integrals,
which appeared in
(\ref{3.12}), and obtain the final result:
\be
\lb{3.13}
Z(N) =  (p+1)^{N} \, \sum_{m=0}^{p} t_{m}^{N} =
(p+1)^{N} \, ( t_{0}^{N} + p \, t_{1}^{N} ) =
 \left( (e^{K} + p)^{N} + p \, (e^{K} -1)^{N} \right) \; .
\ee
It is possible to check this result provided that matrix $V$
has one eigenvalue $(x+p)$ and $p$ eigenvalues $(x-1)$.
Let us note that cyclic symmetry of
(\ref{3.12}) with respect to the permutation of indices
$(1 \rightarrow 2 \rightarrow \dots \rightarrow N \rightarrow 1)$
(evident from representation (\ref{3.3}))
becomes manifest if we put
$a^{1}_{i} = 1$ and $a^{2}_{i} = -1$. Thus, the algebra
$\Pi_{p+1}(2N)$ (\ref{1.35})
is fixed completely and defining relations acquire the form:
$$
\t_{i} \, \t_{j} =
q^{\epsilon_{ij}+\delta_{ij}} \, \t_{j} \, \t_{i}  \; , \;\;
\d_{i} \, \d_{j} =
q^{\epsilon_{ij} + \delta_{ij}} \, \d_{j} \, \d_{i} \; , \;\;
\d_{i} \, \t_{j} = q^{-\epsilon_{ij}} \, \t_{j} \, \d_{i}
+ \delta_{ij} \; , \;\;
$$
\be
\lb{3.14}
\overline{\t}_{i} \, \overline{\t}_{j} = q^{\epsilon_{ij}+\delta_{ij}} \,
\overline{\t}_{j} \, \overline{\t}_{i} \;\;  , \;\;
\overline{\d}_{i} \, \overline{\d}_{j} = q^{\epsilon_{ij} + \delta_{ij}} \,
\overline{\d}_{j} \, \overline{\d}_{i} \; , \;\;
\overline{\d}_{i} \, \overline{\t}_{j} = q^{-\epsilon_{ij}} \,
\overline{\t}_{j} \, \overline{\d}_{i}
+ \delta_{ij} \; , \;\;
\ee
$$
\overline{\t}_{i} \, \t_{j} = q^{-\epsilon_{ij}} \,
\t_{j} \, \overline{\t}_{i} \; , \;\;
\overline{\d}_{i} \, \d_{j} = q^{-\epsilon_{ij}} \,
\d_{j} \, \overline{\d}_{i} \; , \;\;
\overline{\d}_{i} \, \t_{j} = q^{\epsilon_{ij}} \,
\t_{j} \, \overline{\d}_{i} \; . \;\;
\overline{\t}_{i} \, \d_{j} = q^{\epsilon_{ij}} \,
\d_{j} \, \overline{\t}_{i} \; , \;\;
$$
where
$\epsilon_{ij} = \; + 1$ for $i > j$ and
$\epsilon_{ij} = \; - 1$ for $i \leq j$.

PGA with many generators
(\ref{3.14}) gives us a possibility to define
(with the help of integral (\ref{1.25})) the convolution product
of two functions (\ref{1.21a})
$f^{(1)}$ and $f^{(2)}$:
$$
f^{(3)}(\t_{1}, \, \overline{\t}_{3}) =
\int \int \, d\t_{2} d\overline{\t}_{2} \,
f^{(1)}(\t_{1}, \, \overline{\t}_{2}) \,
f^{(2)}(\t_{2}, \, \overline{\t}_{3}) \,
\mu(\t_{2}, \, \overline{\t}_{2}) \; .
$$
This convolution product is equivalent to the matrix
product of the corresponding matrices of the coefficients
$f^{(3)}_{ij} =  f^{(1)}_{ik} \, f^{(2)}_{kj}$.
An analogous convolution product of the function
(\ref{1.21a}) with holomorphic functions
(\ref{1.22}) is equivalent to the action
of a matrix on a vector or covector, while integral
(\ref{1.29}) is nothing but the contraction of the vector and
covector.

Let us introduce the PG coherent states
$< \t|$ and $| \overline{\t} >$,
which can be defined by the equations
\be
\lb{4.2}
 \hat{\d}^{k} \, |\overline{\t}_i>  = |\overline{\t}_i> \,
\overline{\t}^{k}_{i} \; , \;\;
<\t_i | \, \hat{\t}^{k}  = \t^{k}_{i} <\t_{i} | \; .
\ee
Here, we use special notation $\hat{\t}$ and $\hat{\d}$
for the PG operators acting on the states $<.|, \; |.>$,
to distinguish them from the PG numbers
$\t_i$, $\overline{\t}_i$ and the corresponding derivatives
$\d_i$, $\overline{\d}_i$.
We postulate the commutation relations
$$
\hat{\t} \, \t_{i} = q^{\xi} \, \t_{i} \, \hat{\t} \; , \;\;
\hat{\d} \, \t_{i} = q^{-\xi} \, \t_{i} \, \hat{\d} \; , \;\;
\hat{\t} \, \overline{\t}_{i} = q^{-\xi}
 \, \overline{\t}_{i} \, \hat{\t} \; , \;\;
\hat{\d} \, \overline{\t}_{i} = q^{\xi}
\, \overline{\t}_{i} \, \hat{\d} \; ,
$$
where $\xi$ is
a constant.
In this case, one can construct coherent states
(\ref{4.2}) explicitly:
$$
|\overline{\t}_{i}>  =
\sum_{k=0}^{p} \,
\frac{\hat{\t}^{k} }{k_{q}!} \, |0> \overline{\t}_{i}^{k} \;\; , \;\;\;
 <\t_{i}|  =
 \sum_{k=0}^{p} \, \t_{i}^{k} \,
<0| \, \frac{ \hat{\d}^{k}}{k_{q}!} \; ,
$$
and check the relations (for all $i$)
\be
\lb{4.7}
1= \sum_{k=0}^{p} \,
\frac{1}{k_{q}!} \, \hat{\t}^{k} \, |0> \,
<0| \, \hat{\d}^{k}  = \int \int \,
d\t_{i} d\overline{\t}_{i} \,
| \overline{\t}_{i} >< \t_{i} | \,
\mu (\t_{i} \overline{\t}_{i}) \; .
\ee

Let us consider the following operator:
$$
H = \sum_{n=0}^{p} \, h_{n} \, \hat{\t}^{n} \,
(g^{-1/2} \, \hat{\d})^{n}
 = \sum_{n=0}^{p} \, h_{n} \, q^{n(1-n)/4} \, \hat{\t}^{n} \,
g^{-n/2} \, \hat{\d}^{n} \; ,
$$
(where $h_{n}$ are real numbers)
as a Hamiltonian of a quantum mechanical system
on the PG space. The hermitian property
$H = H^{\dagger}$ for this Hamiltonian
can be verified by using the conjugation rules (\ref{1.13}).
The heat kernel
$$
< \t_{0} | \, U(t) \, | \overline{\t}_{N+1} > =
 < \t_{0} | exp(i \, t \, H) | \overline{\t}_{N+1} > \; ,
$$
is represented as an integral in the PG space
\be
\lb{4.14}
\begin{array}{c}
< \t_{0} | \, U(t) \, | \overline{\t}_{N+1} > =
\int \dots \int \,
d\t_{1} d\overline{\t}_{1} \dots
d\t_{N} d\overline{\t}_{N} \,
\exp_{q} (-\t_{1} \overline{\t}_{1}) \dots
\exp_{q} (-\t_{N} \overline{\t}_{N}) \, \\ \\
< \t_{0} | \, U(\Delta) \, | \overline{\t}_{1} > \,
< \t_{1} | \, U(\Delta) \, | \overline{\t}_{2} > \dots
< \t_{N} | \, U(\Delta) \, | \overline{\t}_{N+1} > \; ,
\end{array}
\ee
where the interval of time $t$ is divided into $N$
equal parts $\Delta = t/N$, and we use eq. (\ref{4.7}).
In the limit $N \rightarrow \infty $ we obtain \cite{FIIK}
\be
\lb{4.15}
< \t_{i} | \, U(\Delta) \, | \overline{\t}_{i+1} > =
\sum_{m=0}^{p} \, \frac{\t_{i}^{m} \; \overline{\t}_{i+1}^{m}}{m_q !} \,
t_{m} \; ,
\ee
where
$$
t_{m}
\simeq \exp \left( - i \, \Delta \, m_{q}! \, \sum_{n=0}^{m} \,
\frac{h_{n} \, q^{n(n+1-2m)/4} }{ (m-n)_{q}! } \right) \; .
$$
Relation (\ref{4.15}) is deduced with the help of
(\ref{4.2}). After substitution of (\ref{4.15}) into
(\ref{4.14}) we obtain the formula which reminds
expression (\ref{3.12}) for the partition function $Z(N)$.
The integral formula (\ref{4.14}) is universal
and naturally appeared in other considerations related to the
applications of PGA.
We hope that,
due to equations of the type
$$\exp_{q}(\t_{2}) \, \exp_{q}(\t_{1}) =
\exp_{q}(\t_{2} + \t_{1}) \;\; (\t_{1}\t_{2} = q\t_{2}\t_{1}) \; ,
$$
and nilpotency of $\t_{i}$ and $\overline{\t}_{i}$,
one can obtain for (\ref{4.14})
the closed expression as a path integral in PG space
in the continuous limit $\Delta \rightarrow 0$
(for every case $p=2,3, \dots$).

At the end of this report, we construct finite dimensional
representations for the quantum group
$GL_{q^{1/2}}(N)$, where $N=2$, by using
the realization of their generators in terms of
elements of PGA. We hope that this construction
can be generalized to the case of arbitrary $N$.
The elements $a,b,c,d$ of
$GL_{q^{1/2}}(2)$ are usually arranged as the $2 \times 2$ matrix
$T_{ij} =
 \left(
\begin{tabular}{c}
$ a, \; b$    \\
$ c, \; d$
\end{tabular}
\right)$,
and commutation relations for these elements are \cite{FRT}
\be
\lb{5.2}
\begin{array}{c}
a \, b = q^{1/2} \, b \, a \; , \;\; a \, c = q^{1/2} \, c \, a \; , \;\;
b \, d = q^{1/2} \, d \, b \; , \;\; c \, d = q^{1/2} \, d \, c \; , \;\; \\
b \, c = c \, b \; , \;\; [a, \, d] = (q^{1/2}-q^{-1/2}) \, b \, c \; .
\end{array}
\ee
The following representation of the
$GL_{q^{1/2}}(2)$ generators in terms of the
$q$-oscillator $\d\t = 1 + q\t\d$ (PG variables) can be obtained
\be
\lb{5.8}
a=g^{\alpha} \, \d \; , \;\; b = \beta \, g^{1/2} \; , \;\;
c= \gamma \, g^{1/2} \; , \;\;
d = \beta \, \gamma \,(q^{1/2}-q^{-1/2})
 \, \t \, g^{-\alpha}  \; ,
\ee
where $\alpha, \beta$ and $\gamma$ are constants.
The quantum determinant (central element for the algebra
$Fun(GL_{q^{1/2}}(2))$) acquires the form
$$
det_{q^{1/2}}(T) = a \, d - q^{1/2} \, b \, c =
-q^{-1/2} \, \beta \, \gamma \; .
$$
The case of the
group $SL_{q^{1/2}}(2)$ corresponds to the choice
$det_{q^{1/2}}(T)=1$ and leads to an additional constraint
on the parameters $\beta \, \gamma = - q^{1/2}$.
Now one can use the finite dimensional representations
(\ref{1.9a}) for PGA and construct (with the help
of (\ref{5.8})) the $(p+1)$- dimensional representations
for the algebras
$Fun(GL_{q^{1/2}}(2))$ and $Fun(SL_{q^{1/2}}(2))$ when
$q^{p+1}=1$.

We note that formulas (\ref{5.8}) give also the
realization of the groups
$GL_{q^{1/2}}(2)$ and $SL_{q^{1/2}}(2)$ for arbitrary parameters $q$,
if we do not require the nilpotency conditions for
$\t$ and $\d$.
The infinite dimensional representations of such $q$-oscillator
define infinite dimensional representations of
$GL_{q^{1/2}}(2)$ and $SL_{q^{1/2}}(2)$. These representations
are related to those considered in \cite{VS}.

\end{document}